\documentclass[10pt]{article}
\usepackage{amssymb, amsmath}
\usepackage{graphicx}
\usepackage[latin1]{inputenc}
\usepackage{natbib}
\usepackage{epsfig, epic, eepic}

\newtheorem{theorem}{Theorem}

\newtheorem{lemma}{Lemma}

\newtheorem{notation}{Notation}
\numberwithin{equation}{section}

\begin{document}


\newcommand{\calA}{{\cal A}}
\newcommand{\calB}{{\cal B}}
\newcommand{\calF}{{\cal F}}
\newcommand{\calG}{{\cal G}}
\newcommand{\calR}{{\cal R}}
\newcommand{\calZ}{{\cal Z}}
\def\D{\mathcal{D}}
\def\L{\mathcal{L}}
\def\S{\mathcal{S}}
\def\I{\mathcal{I}}
\def\V{\mathcal{V}}
\def\E{\mathcal{E}}

\newcommand{\N}{{\mathbb N}}
\newcommand{\Z}{{\mathbb Z}}
\newcommand{\Q}{{\mathbb Q}}
\newcommand{\R}{{\mathbb R}}
\newcommand{\C}{{\mathbb C}}
\newcommand{\K}{{\mathbb K}}
\newcommand{\kk}{{\mathrm k}}

\def\J{\mathrm{J}}
\def\x{\mathrm{x}}
\def\a{\mathrm{a}}
\def\d{\mathrm{d}}

\def\GL{\mathrm{GL}}
\def\det{\mathrm{det}}
\def\SL{\mathrm{SL}}
\def\PSL{\mathrm{PSL}}
\def\PGL{\mathrm{PGL}}
\def\O{\mathrm{O}}

\def\gl{\mathfrak{gl}}
\def\g{\mathfrak{g}}
\def\h{\mathfrak{h}}

\newcommand{\Frac}[2]{\displaystyle \frac{#1}{#2}}
\newcommand{\Sum}[2]{\displaystyle{\sum_{#1}^{#2}}}
\newcommand{\Prod}[2]{\displaystyle{\prod_{#1}^{#2}}}
\newcommand{\Int}[2]{\displaystyle{\int_{#1}^{#2}}}
\newcommand{\Lim}[1]{\displaystyle{\lim_{#1}\ }}

\newenvironment{menumerate}{%
    \renewcommand{\theenumi}{\roman{enumi}}%
    \renewcommand{\labelenumi}{\rm(\theenumi)}%
    \begin{enumerate}} {\end{enumerate}}
     
\newenvironment{system}[1][]%
	{\begin{eqnarray} #1 \left\{ \begin{array}{lll}}%
	{\end{array} \right. \end{eqnarray}}

\newenvironment{meqnarray}%
	{\begin{eqnarray}  \begin{array}{rcl}}%
	{\end{array}  \end{eqnarray}}

\newenvironment{marray}%
	{\\ \begin{tabular}{ll}}
	{\end{tabular}\\}

\newenvironment{program}[1]%
	{\begin{center} \hrulefill \quad {\sf #1} \quad \hrulefill \\[8pt]
		\begin{minipage}{0.90\linewidth}}
	{\end{minipage} \end{center} \hrule \vspace{2pt} \hrule}

\newcommand{\entrylabel}[1]{\mbox{\textsf{#1:}}\hfil}
\newenvironment{entry}
   {\begin{list}{}%
   	{\renewcommand{\makelabel}{\entrylabel}%
   	  \setlength{\labelwidth}{40pt}%
   	  \setlength{\leftmargin}{\labelwidth + \labelsep}%
   	}%
   }%
   {\end{list}}

\newenvironment{example}{\par \noindent {\bf Example. }}
			{\hfill $\blacksquare$ \par}

\newenvironment{Pmatrix}
        {$ \left( \!\! \begin{array}{rr} } 
        {\end{array} \!\! \right) $}

\newcommand{\fleche}[1]{\stackrel{#1}\longrightarrow}
\def\ssi{si et seulement si\ }
\newcommand{\tab}{\hspace*{\fill}}
\newcommand{\bs}{{\backslash}}
\newcommand{\eps}{{\varepsilon}}
\newcommand{\into}{{\;\rightarrow\;}}
\newcommand{\PD}[2]{\frac{\partial #1}{\partial #2}}
\def\Hat{\widehat}
\def\Bar{\overline}
\def\vect{\vec}
\def\fbar{{\bar f}}
\def\xbar{{\bar \x}}
\newcommand{\afaire}[1]{$$\vdots$$ \begin{center} {\sc #1} \end{center} $$\vdots$$ }
\newcommand{\pref}[1]{(\ref{#1})}

\def\Maple{{\sc Maple}}
\def\RG{{\sc Rosenfeld-Gr\"obner}}



\newcommand{\algf}{\sffamily}
\newcommand{\BEGIN}{{\algf begin}}
\newcommand{\END}{{\algf end}}
\newcommand{\IF}{{\algf if}}
\newcommand{\THEN}{{\algf then}}
\newcommand{\ELSE}{{\algf else}}
\newcommand{\ELIF}{{\algf elif}}
\newcommand{\FI}{{\algf fi}}
\newcommand{\WHILE}{{\algf while}}
\newcommand{\FOR}{{\algf for}}
\newcommand{\DO}{{\algf do}}
\newcommand{\OD}{{\algf od}}
\newcommand{\RETURN}{{\algf return}}
\newcommand{\PROCEDURE}{{\algf procedure}}
\newcommand{\FUNCTION}{{\algf function}}
\newcommand{\INDENTER}{{\algf si} \=\+\kill}

\newcommand{\target}{\mathop{\mathrm{t}}}
\newcommand{\source}{\mathop{\mathrm{s}}}
\newcommand{\trdeg}{\mathop{\mathrm{tr~deg}}}
\newcommand{\jet}[2]{\jmath_{#1}^{#2}}
\newcommand{\rank}{\operatorname{rank}}
\newcommand{\sign}{\operatorname{sign}}
\newcommand{\ord}{\operatorname{ord}}
\newcommand{\aut}{\operatorname{aut}}
\newcommand{\Hom}{\operatorname{Hom}}
\newcommand{\codim}{\operatorname{codim}}
\newcommand{\coker}{\operatorname{coker}}
\newcommand{\rp}{\operatorname{rp}}
\newcommand{\leader}{\operatorname{ld}}
\newcommand{\card}{\operatorname{card}}
\newcommand{\Fr}{\operatorname{Frac}}
\newcommand{\RF}{\operatorname{\mathsf{reduced\_form}}}
\newcommand{\rang}{\operatorname{rang}}

\def \Id{\mathrm{Id}}

\def \diff{\mathrm{Diff}^{\mathrm{loc}} }
\def \diffg{\mathrm{Diff} }
\def \Esc{\mathrm{Esc}}

\newcommand{\initial}{\mathop{\mathsf{init}}}
\newcommand{\separant}{\mathop{\mathsf{sep}}}
\newcommand{\quo}{\mathop{\mathsf{quo}}}
\newcommand{\pquo}{\mathop{\mathsf{pquo}}}
\newcommand{\lcoeff}{\mathop{\mathsf{lcoeff}}}
\newcommand{\mvar}{\mathop{\mathsf{mvar}}}

\newcommand{\prem}{\mathop{\mathsf{prem}}}
\newcommand{\remp}{\mathrel{\mathsf{partial\_rem}}}
\newcommand{\remf}{\mathrel{\mathsf{full\_rem}}}
\renewcommand{\gcd}{\mathop{\mathrm{gcd}}}
\newcommand{\pairs}{\mathop{\mathrm{pairs}}}
\newcommand{\dd}{\mathrm{d}}
\newcommand{\ideal}[1]{(#1)}
\newcommand{\cont}{\mathop{\mathrm{cont}}}
\newcommand{\pp}{\mathop{\mathrm{pp}}}
\newcommand{\pgcd}{\mathop{\mathrm{pgcd}}}
\newcommand{\ppmc}{\mathop{\mathrm{ppcm}}}
\newcommand{\init}{\mathop{\mathrm{initial}}}

\bibliographystyle{amsalpha}

\title{On the geometry of the first and second Painlevé equations}%
\author{
	Raouf Dridi \\
{\small 
{ {\sf 
	dridi.raouf@gmail.com
}}}
\\
{\small 
 {\it 
	Doppler Institute for Mathematical Physics and Applied Mathematics, 
 }
}
\\
{\small 
 {\it 
	Brehova 7, 11519 Prague 1, Czech Republic
 }
}
}

\date{\today}
\maketitle

\begin{abstract}
	In this  paper 
	we \emph{explicitly} compute the transformation that maps the generic second order differential
	equation~$y''= f(x, y, y')$ to the Painlevé first equation~$y''=6y^2+x$   (resp.  the Painlevé second
	equation ${y''=2 y^{3}+yx+ \alpha}$).   This change of coordinates, which is function of~$f$ and its
	partial derivatives, does not exist for every~$f$; it is necessary that the function $f$ satisfies certain conditions that define
	the equivalence class of the considered Painlevé equation.
	In this work we won't consider these conditions and the existence issue is solved  \emph{on line} as follows:  
	If  the input equation is known then it suffices to specialize the change of coordinates
	on this equation and test by simple substitution if the equivalence holds. 
	The other innovation of this work lies in the exploitation of discrete symmetries for solving the equivalence problem.
\end{abstract}

\section{Introduction}\label{}
By fiber-preserving transformations we mean analytical transformations of the form   
$$
	\C^2 \ni (x,y)\into \left(\bar x(x),\ \bar y(x,y) \right) \in \C^2
$$
with the condition $\bar x_x\bar y_y\neq0$  expressing their  local invertibility. These 
transformations form a Lie pseudo-group with 
\begin{equation}\label{lie}
	\bar x_y = 0,\, \bar x_x\bar y_y\neq0
\end{equation}
as defining system.

As indicated in the abstract, our aim is to \emph{explicitly} compute the  transformation of this form that
 maps the second order equation 
$$
 	\E_f: \, y''=f(x,y,y'), 
$$
where $y' = \Frac{d}{dx} y(x)$,  to  the first Painlev\'e equation
 (resp. to the second Painlevé equation). This change of coordinates, which is clearly function of $f$ and its
partial derivatives, does not exist for every~$f$; it is necessary that the function $f$ satisfies certain conditions that define
the equivalence class of the considered Painlevé equation. Comparing to   \cite{kamran85} and \cite{Hietarinta02}, 
the existence issue is solved here
 \emph{on line}  as follows:  
If  the input equation is known then it suffices to specialize the change of coordinates
 on this equation and test by simple substitution if the equivalence holds.

The calculations of transformation   \emph{candidates}  are based on the following result~\cite{dridi:issac07}. 
 Given a Lie pseudo-group of  transformations $\Phi$ and denote by $\S_{\E_f, \Phi}$ 
the symmetry pseudo-group  of the equation~$\E_f$ w.r.t to $\Phi$ i.e., $\S_{\E_f, \Phi}=\Phi\cap\diff(\E_f)$. In~\cite{dridi:issac07}, we proved

 \begin{menumerate}

\item	[$(i)$]
The number of  constants appearing in the change of coordinates  is exactly the dimension of $\S_{\E_\fbar, \Phi}$.  
This implies that when this dimension vanishes the change of coordinate can be obtained without
  integrating differential equations.  Also, we have $\dim(\S_{\E_f, \Phi})=\dim(\S_{\E_\fbar, \Phi})$.

\item	[$(ii)$]
In the particular case when $\dim(\S_{\E_\fbar, \Phi})=0$,  the transformation $\varphi$ is algebraic in $f$
and its partial derivatives and it is obtained \emph{without} solving differential equations.
The degree of this transformation $\varphi$ is exactly equal to the finite value $\card(\S_{\E_\fbar, \Phi})$. 
\end{menumerate}

The last case is exactly what happens when $\E_\fbar$ is one of the Painlevé equations and~$\Phi$ is the
pseudo-group of fiber-preserving transformations or more generally point transformations.
Indeed,  the classical Lie analysis shows that the  point symmetry pseudo-group of each one
 of Painlevé equations is zero-dimensional.  Moreover, according to the fact that 
 the unique transformations that preserve the  singularity structure are  the homographic transformations, one can show by
 straightforward computations that the point symmetry pseudo-group of Painlevé one is 
 \begin{system}\label{sym_gr_p1}
 \bar x &=& x\Frac{\bar y^2}{y ^2} ,  \\[3mm]
 \bar y^ 5&=& y ^5,
 \end{system}%
%
and  
 \begin{system}\label{sym_gr_p2}
 \bar x &=& x\Frac{\bar y^2}{y ^2},  \\[4mm]
 \bar y^ 3&=& \Frac{ \bar \alpha }{\alpha}  y ^3,\\[3mm]
 \bar \alpha^2 &=& \alpha^2,
 \end{system}%
for  Painlevé two when $\alpha\neq 0$ and
 \begin{system}\label{sym_gr_p2_alpha_zero}
 \bar x &=& x\Frac{\bar y^2}{y ^2},  \\[3mm]
 \bar y^ 6&=&  y ^6,
 \end{system}%
when $\alpha=0$.



Fiber-preserving transformations are suited when dealing with Painlevé equations. In particular, such transformations
preserve the integrability in the sense of Poincaré \cite{conte}.
However, since Painlevé equations lie in the class of equations of the form
$$
 	y''= A(x, y)  + B(x, y)y' + C (x, y) y'^2  + D(x, y)y'^3
$$
which is invariant under point transformations\footnote{Indeed, as remarked by E. Cartan \cite{cartan24}, the above equation
can  always be regarded as the geodesics equation of a projective structure on a surface
with local coordinates $x$ and $y$ and thus invariant under point transformations.}, 
we  consider in the last section of this paper the equivalence under these more general transformations.

\section{Building the invariants}
Le $(x,\,y,\,p=y')$ be a system of local coordinates of  $\J^1 = \J^1(\C,\C)$,
the space of first order jets of functions $\C\ni x\into y(x)\in \C$~\cite{olver}. Two scalar second order ordinary equations 
 $$
	    \E_f :  y''=f(x,\,y,\,y')  \mbox{ and } \E_\fbar :  
	    \bar y'' = \bar f(\bar x,\,\bar y,\, \bar y')
 $$
 are said to be equivalent under a point transformation  $\varphi$  if its first prolongation~(to~$\J^1$) maps the contact forms
 \begin{system}\nonumber
\omega^1 &=& dy-p\ dx \\
\omega^2 &=&  dp-f(x,y,p)\ dx
 \end{system}%
to  the contact forms
 \begin{system}\nonumber
\bar \omega^1 &=& d\bar y-\bar p\ \bar dx \\
\bar \omega^2 &=&  d\bar p-\bar f(\bar x,\bar y,\bar p)\ d\bar x
 \end{system} %
 up to an invertible $2\times 2$-matrix of the form
 $$
         \begin{pmatrix}
          a_1 & 0 \\
          a_2 & a_3  \\
        \end{pmatrix}.
 $$
The $a_i$ are functions  from $\J^1$ to $\C$.
 To encode equivalence under fiber-preserving transformations~(i.e. taking in account the Lie equations \pref{lie}) we must have $$\varphi^* d\bar x = a_4 dx$$  for a certain function $a_4 : \J^1\into\C$. Summarizing: 
two second order differential equations~$\E_f$ and~$\E_\fbar$ are equivalent under a fiber-preserving transformation $\varphi$ if and only if
\begin{equation*}\label{pff}
\varphi^*
        \begin{pmatrix}
          d\bar{y}-\bar{p}\ d\bar{x} \\
          d\bar{p}-\bar{f}(\bar x,\bar y,\bar p)\ d\bar{x} \\
          d\bar{x}
        \end{pmatrix}%
        =
          \begin{pmatrix}
          a_1 &0 & 0 \\
          a_2 & a_3 & 0 \\
          0  & 0& a_4
        \end{pmatrix}%
        \begin{pmatrix}
          dy-p\ dx \\
          dp-f(x,y,p)\ dx \\
          dx
        \end{pmatrix}.%
\end{equation*} 
For this problem,  Cartan's equivalence method \cite{olver95} gives  three fundamental invariants 
\begin{system}\nonumber
I_{3}=-\Frac{f_{ppp} a_{4}}{2{a_{1}}^{2}}, \\  [4mm]
I_{2}= \Frac{f_{yp} - D_{x}f_{pp}}{2a_{1} a_{4}},\\ [4mm]
I_{1}= \Frac{(2 f_{yy} - D_{x}f_{yp} - f_{pp} f_{y} + f_{yp} f_{p}) a_{1} + ( - f_{yp} + D_{x}f_{pp}) a_{4} a_{2}}{2{a_{1}}^{2} {a_{4}}^{2}}
\end{system}%
and six invariant derivations defined on certain manifold $\tilde M$, fibred over~$\J^1$, with local coordinates de~${(x,y,p,a_1,a_2,a_4)}$. Here,~
$D_x={\partial_x} + p{\partial_y} +f{\partial_p}$ is Cartan vector field.  

When specializing on the Painlev\'e equations,
 the two fundamental invariants~$I_2$ and~$I_3$ vanish. On this splitting   branch, the application of
  Jaccobi identity  to the final structure equations shows that  among the six invariant derivations only two can produce new invariants. These two derivations are
\begin{system} \nonumber
	X_1 & = & \Frac{1}{a_{1}}\partial_y - \Frac{a_{2} a_{4}}{{a_{1}}^{2}}\partial_p 
	- \Frac{1}{2}f_{pp} \partial_{a_{1}} -  \Frac{1}{2}\Frac{f_{py} }{a_{4}}\partial_{a_{2}}, \\ [5mm]
	X_3 & = & \Frac{1}{a_{4}}\partial_x + \Frac{p}{a_{4}}\partial_y + 
	\Frac{f }{a_{4}}\partial_p + a_{2}\partial_{a_{1}} - 
	\Frac{f_{y}  a_{1}}{{a_{4}}^{2}}\partial_{a_{2}} 
	+ \frac{2 a_{2} a_{4} + f_{p}  a_{1}}{a_{1}}\partial_{a_{4}}. \\ 
\end{system}%

\begin{notation}
In the sequel, $I_{1;j\cdots k}$  denotes the differential invariant $X_k\cdots X_j(I_1)$.  For instance, the invariant $I_{1;33}$ is obtained by differentiating twice the fundamental invariant $I_1$ with respect to invariant derivation $X_3$.
\end{notation}

\section{The first Painlev\'e equation $ y'' = 6 y^2+  x$}  

Since the associated fiber-preserving   symmetry Lie pseudo-group is zero dimensional, this justify the following lemma:
\begin{lemma}
	The specialization of the invariants
$$ 
	I_1,\ I_{1;3},\ I_{1;33},\ 
    	\frac{I_{1;333}}{I_{1;33}},\
    	{\frac { I_{{1;3333}}}{ I_{{1;33}}}}-{\frac {43}{120}} I_{{1;33}},\
    	{\frac { I_{1;33333}}{ I_{{1;33}}}}-\frac{5}{4} I_{{1;33}}
$$
on the first Painlev\'e equation gives six invariants functionally independent defined on~$\tilde M$.
\end{lemma} 
	The problem with the above invariants is that they do depend on extra parameters $ a_1$, $ a_2$ and~$ a_4$. Fortunately, in our
	zero-dimensional case, we can normalize
	(e.i. eliminate) these parameters by setting
\begin{equation}\label{a}
 	 I_1=-12,\  I_{1;3}= 0, \frac{I_{1;333}}{ I_{1;33}}=1.
 \end{equation}
 	Now substituting the values of the parameters in the remaining invariants give us, due again to our zero-dimensional case,
	three  functionally independent invariants  now do not depending on the extra parameters.
	
	
	Writing  the equality of the invariants and simplifying  the obtained system, by computing a characteristic set  \cite{kolchin:livre, BLOP},  gives an algebraic transformation of degree~5 
\begin{system}\label{cara}
\bar p  &=& 129600 {\Frac {\left( 5 {I_{1;33}}^{2}+4 I_{1;33333} \right) }{{I_{1;33}}^{3}}} {\bar y }^{4}\\  
\bar x  &=&  -6 \Frac { \left( 120 I_{1;3333}+43 {I_{1;33}}^{2} \right) }{I_{1;33}^2}\bar y^2,\\ 
{\bar y}^{5} &=& -{\Frac {1}{23328000}} {\Frac { {I_{1;33}}^
{5}}{
\left( 5 {I_{1;33}}^{2}+4 I_{1;33333} \right)^2
}}. 
\end{system}%
	In these formulae the invariants are normalized using  (\ref{a}), that is, do not depending on the extra parameters.
	 According to  (ii) of the introduction and \pref{sym_gr_p1},  we have

\begin{theorem}\label{prop1}
A second order differential equation  $\E_f$  is equivalent to the first Painlev\'e equation by a fiber-preserving transformation if and only if this transformation is given by (\ref{cara}) and the normalization~(\ref{a}).
\end{theorem}

Let us explain how Theorem 1 can be used in practice.  Consider the following equations
\begin{equation}\label{ref_eq1}
	y'' = c \frac{y'^2}{y}  + \frac{1}{y} (y^4+x), 
\end{equation}
and
\begin{equation}\label{ref_eq2}
	y'' = c \frac{y'^2}{y}  + y (y^4+x) .
\end{equation}
The question is to determine the values of the parameter $c$ for which the above equations
can be mapped to the first Painlevé equation (and compute the equivalence transformation
when the equivalence holds).  

First of all, the fact that the derived invariants $I_{1; 1}$ vanishes on the  first Painlevé equation  restricts
the possible values of $c$ to $\{-1,\, 3\}$ for the first equation and to $\{-3,\, 5\}$ for the second equation.

The second  step is to specialize \pref{cara} on the given equation to obtain transformation candidates.
In step 3,   we have to check whether the pullback of the first Painlevé equation w.r.t these
candidates is exactly the considered  equation.

In the case of  equation \pref{ref_eq1}, the specialization yields
\begin{system}\label{tr_eq1_c_1}
\bar p &=& 36\,{\Frac {{\bar y}^{4}p}{{y}^{7}}},\\[5mm]
\bar x  	  &=&6\,{\Frac {{\bar y}^{2}x}{{y}^{4}}},\\[5mm]
{\bar y}^{5} &=&{\Frac {1}{108}}\,{y}^{10}
\end{system}%
for $c=-1$ and 
\begin{system}\label{tr_ref_eq_c_3}
\bar p &=&-864\,{\Frac {{\bar y}^{4}{y}^{5} \left( 625\,{x}^{5}-2079 \right) 
 \left( -25\,y{x}^{3}+250\,p{x}^{4}+21\,{y}^{3} \right) }{ \left( 50\,
{x}^{3}+3\,{y}^{2} \right) ^{4}}},\\[5mm]
\bar x &=& 6\,{\Frac { \left( 2500\,{x}^{5}-891 \right) {y}^{4}{\bar y}^{2}}{
 \left( 50\,{x}^{3}+3\,{y}^{2} \right) ^{2}}},\\[5mm]
{\bar y}^{5} &=& -{\Frac {1}{ 31104}}\,{\Frac { \left( 50\,{x}^{3}+3\,{y}^{2} \right) ^{5}}{{y}^{10}
 \left( 625\,{x}^{5}-2079 \right) ^{2}}}
\end{system}%
for $c=3$.    The  third steps shows that the equivalence holds only for $c=-1$
 and  the equivalence transformation is \pref{tr_eq1_c_1}.  We can also deduce, according
to $(ii)$ in the introduction, that the cardinal of the fiber-preserving (point)
symmetry group  of  the equation \pref{ref_eq1} with $c=-1$  is equal to 10. 

The same calculations show that equation \pref{ref_eq2} can not be mapped to the first 
Painlevé equation. In particular, we have a {\sf division by zero error} in step 2 for $c=5$. {\it Warning:}
This error doesn't mean that the method failed. In fact it is part of the method and implies that
no equivalence transformation does exist.

Time estimates are given in tables bellow where $P_1$ refers to the first Painlevé equation.

{\small
\begin{center}
\begin{tabular}{|c|c|c|}
\hline
 & { Computation of transformation candidates} &  { Checking equivalence with $P_1$}  \\
 \hline
$c = -1$ & 0.15 & (yes)  0.04  \\
 \hline
$c = 3$ & 2.13 & (no) 0.13  \\
\hline
\end{tabular}
\label{table1} \\
{\bf Table 1. } {\it Time estimates (in seconds) for  $y'' = c \Frac{y'^2}{y}  + \Frac{1}{y} (y^4+x) $}
\end{center}

\begin{center}
\begin{tabular}{|c|c|c|}
\hline
 & Computation of transformation candidates &  Checking equivalence with $P_1$  \\
 \hline
$c = -3$ & 0.35 & (no)  0.03  \\
 \hline
$c = 5$ & {\sf division by zero error} & (no) 0.00   \\
\hline
\end{tabular}
\label{table2} \\
{\bf Table 2. }  {\it Time estimates (in seconds) for  $y'' = c \Frac{y'^2}{y}  +  y(y^4+x) $}
\end{center}
}


\section{The second Painlev\'e equation $ y''=2 y^{3}+  xy + \alpha$}

Again, due to the zero-dimensionality,  there exists seven invariants defined on the manifold of local coordinates
 $(x,y,p, a_1, a_2,a_4,\alpha)$ such that when specialized, on Painlevé two, they give exactly seven functionally independent functions. For instance, one can take the invariants
  $I_1, I_{1;3}, I_{1;31}, I_{1;33}, I_{1;331}, I_{1;3331}$ and $I_{1;33311}$.  We normalize   $ a_1,  a_2$ and $ a_4$ by setting
\begin{equation}\label{normalIP2fiber}
I_1 =-12, \  
I_{1;3} = -12, \
I_{1;31} =  0,
\end{equation}
and as in the previous section, we obtain 
 \begin{system}\label{cara1p2}
 \bar p &=& \Frac{1}{6} \left( {\Frac { I_{1;33311} \left( I_{1;3331} +4032 \right) }{I_{1;33311}I_{1;33} -3096576
-4032 I_{1;331}}} \right){\bar y }^{2}\bar \alpha,\\[5mm]
\bar x &=& - \left ( 16 +{\Frac {1}{72}} I_{1;331} \right ){\bar y }^{2},\\
{\bar y }^{3} &=& 48384 {\Frac {\bar \alpha }{I_{1;33311}I_{1;33} -3096576-4032 I_{1;331}}},\\
{\bar \alpha }^{2} &=&- {\Frac{1} {112 I_{1;33311}  
 \left( 16257024+8064  I_{1;3331}  +{ I_{1;3331}  }^{2}
 \right) }}
(
{ I_{1;33311}  }^{2}{ I_{1;33}  }^{2}\\[5mm]
&& -8064  I_{1;33311}   I_{1;33}
I_{1;331} -6193152  I_{1;33311}   I_{1;33}  \\[3mm]
&& +9588782923776+
24970788864  I_{1;331}  +16257024 { I_{1;331}  }^{2}) . 
 \end{system}%
 when $\alpha\neq 0$ and
\begin{system}\label{cara1p2_alpha_zero}
 \bar p  &=&\Frac {1}{{290304}}\,{y }^{5}{I_{1;33311}} \left( 4032
+I_{1; 3331} \right) ,\\[3mm]
x  &=& -{\Frac {1}{72}}\,{y}^{2}
 \left( 1152+I_{1;331}\right) ,\\
 {y }^{6} &=&-20901888\,{\Frac {1}{I_{1;33311} \left( 4032+I_{1;3331} \right) ^{2}}}, 
\end{system}%
 when $\alpha = 0$. The comparison with  the symmetry pseudo-groups \pref{sym_gr_p2} and \pref{sym_gr_p2_alpha_zero}  proves
\begin{theorem}
A second order differential equation can be mapped to the second Painlevé equation  $ y''=2 y^{3}+yx+ \alpha$  by a fiber-preserving transformation if and only if this transformation is given by~\pref{cara1p2} if $\alpha\neq 0$ and by   \pref{cara1p2_alpha_zero} otherwise with the normalization \pref{normalIP2fiber} in both cases.
\end{theorem}
Let us remark that  \pref{cara1p2_alpha_zero} can be obtained from  \pref{cara1p2} (as well as \pref{sym_gr_p2_alpha_zero} from \pref{sym_gr_p2})
 by eliminating the $\bar \alpha$ and taking in account the  functional dependence between the invariants resulting from $\bar \alpha=0$.
Nevertheless,  it is more safe  to separate the two cases ($\alpha \neq 0$ and $\alpha=0$).
 
\section{Equivalence under point transformation}
The equivalence  problem under the more general point transformations
 naturally arises since Painlevé equations  belong to the class of equations of the form
$$
 	y''= A(x, y)  + B(x, y)y' + C (x, y) y'^2  + D(x, y)y'^3
$$
which is invariant under point transformations. 
In this case our starting Pfaffian system is
$$
\varphi^*
        \begin{pmatrix}
          d\bar{y}-\bar{p}\ d\bar{x} \\
          d\bar{p}-\bar{f}(\bar x,\bar y,\bar p)\ d\bar{x} \\
          d\bar{x}
        \end{pmatrix}%
        =	
        \begin{pmatrix}
          a_1 &0 & 0 \\
          a_2 & a_3 & 0 \\
          a_4  & 0& a_5
        \end{pmatrix}%
        \begin{pmatrix}
                  dy-p\ dx \\
          dp-f(x,y,p)\ dx \\
          dx
        \end{pmatrix}%
$$
for which we  normalize $a_3$ and prolong to obtain involution and four fundamental invariants defined on 8-dimensional manifold. For the above class, only two invariants are not identically zero 
\begin{eqnarray}\nonumber
K_1 &=& (6 f_{yy} - 4 D_{x}f_{yp} + {D_{x}}^{2}f_{pp} - 3 f_{y} f_{pp} + 4 f_{yp} f_{p} - D_{x}f_{pp} f_{p})/(a_{1} {a_{5}}^{2}),\\\nonumber
K_2 &=& (2 f_{y} f_{ppp} a_{5} + 4 f_{yp} f_{p} a_{4} - D_x f_{pp} f_{p} a_{4} - 3 f_{y} f_{pp} a_{4} - a_{5} f_{pp} f_{yp} + a_{5} f_{pp} D_x f_{pp} + 6 a_{4} f_{yy} \\\nonumber
&& + a_{4} D_x D_x f_{pp}  - a_{5} D_x f_{ppp} f_{p} - a_{5} f_{ppp} D_x f_{p} - 4 a_{4} D_x f_{yp} - 2 f_{yyp} a_{5} + 2 a_{5} D_x f_{ypp}\\\nonumber
&& - a_{5} D_x D_x f_{ppp})/({a_{5}}^{2} {a_{1}}^{2}).
\end{eqnarray}
 As  in the fiber-preserving case, only two invariant derivations $X_1$ and $X_3$ (one page long) are needed.

\begin{theorem}
A second order ordinary differential equation $y'' = f(x,\, y,\, y')$ is equivalent  

(i) to the first Painlevé equation $y''= 6y^2 +x$ under a point transformation if and only if this transformation is given by
\begin{system}
\bar p &=&{\Frac {5}{1056}}\,{\Frac { \left( 2^{15} 3^5 11^3 \,K_{{1;33333}}+{K_{{1;33313}}}^{3} \right) }{{K_{{1;33313}}}^{2}}}{\bar y}^{4},\\[5mm]
\bar x &=& -6\,{\Frac { \left( 2^9 3^3 5 11^2\,K_{{1;3333}}
+43\,{K_{{1;33313}}}^{2} \right) }{{K_{{1;33313}}}^{2}}}{\bar y}^{2},\\[5mm]
{y}^{5} &=&-{\Frac {88}{375}}\,
{\Frac {{K_{{1;33313}}}^{5}}{
\left( 2^{15} 3^5 11^3 \,K_{{1;33333}}+{K_{{1;33313}}}^{3} \right)^2
}}
\end{system}%
with the normalization 
$$
K_1=-12,\, K_2=0,\, K_{1;1}=0,\, K_{1;3}=0,\, K_{1; 33}/K_{1;333}=720.
$$

(ii)  to the second Painlevé equation $ y''=2 y^{3}+  xy + \alpha$ under a point transformation if and only if this transformation is given by
\begin{system}\label{p2_pt}
\bar p  &=&-{\Frac {1}{18}}{\Frac { K_{2;3} \left( 15\,K_{2;3}K_{1;33}-216000+4032\,K_{2;3}-450\,K_{1;331}-
50\,K_{2;3}K_{1;333} \right) }{25\,K_{2;3}K_{1;33}-115200+1728\,K_{2;3}-150\,K_{1;331}}} {\bar y }^{2}\bar \alpha,\\[3mm]
\bar x &=&{\Frac {1}{3600}}  \left( 25\,K_{2;3}K_{1;33}+336\,K_{2;3}-57600-50\,K_{1;331} \right){\bar y_{{}}}^{2} ,\\[3mm]
 {\bar y }^{3} &=&-1800\,{\Frac {\bar \alpha }{25\,K_{2;3}K_{1;33}-115200+1728\,K_{2;3}-150\,K_{1;331}}},\\
{\bar \alpha }^{2}&=&-108\,{\Frac { \left( 25\,K_{2;3}K_{1;33}-
115200+1728\,K_{2;3}-150\,K_{1;331} \right) ^{2}}{K_{2;3}
 \left( 15\,K_{2;3}K_{1;33}-216000+4032\,K_{2;3}-450\,K_{1;331}-50\,K_{2;3}K_{1;333} \right) ^{2}}}
\end{system}%
when $\alpha\neq 0$ and 
\begin{system}\label{p2_pt_alpha_zero}
\bar p  &=& -{\Frac {1}{16200}}K_{2;3  } \left( 576\,K_{2;3}+25\,K_{2;3  }K_{1;333 }+30\,K_{2;3  }K_{1;33}-64800
 \right)\,{\bar y }^{5},\\[3mm]
 \bar x  &=& {\Frac {1}{1080}}\left( 5\,K_{2;3  }K_{1;33  }-5760-72\,K_{2;3  } \right)\,{\bar y }^{2} ,\\[3mm]
{\bar y }^{6} &=& - {\Frac {87480000}{K_{2;3  } \left( 576\,K_{2;3  }+25\,K_{2;3  }K_{1;333}+
30\,K_{2;3  }K_{1;33}-64800 \right) ^{2}}}
\end{system}%
when $\alpha = 0$,
with the normalization 
$$
K_1=-12,\, K_2=0,\, K_{1;1}=0,\, K_{1;3}=0,\, K_{2;3}/K_{1;31}=-5/24.
$$
\end{theorem}

\begin{example}
Let us terminate with considering the equivalence of the two equations \pref{ref_eq1} and
\pref{ref_eq2} with the second Painlevé equation under point transformations. Here, computations 
are done with arbitrary~$c$.  

{\it The equation \pref{ref_eq1}:}  
The specialization of \pref{p2_pt_alpha_zero} on this equation yields (after 0.512 seconds) a transformation candidate depending on~$c$
and which is too long to include in this paper.  The variable~$\bar x$  doesn't depend on $p$ in only
tow cases $c\in \{-1, 3\}$ and this two values return a {\sf division by zero error}
when computing the others components. The same thing happens with the specialization of \pref{p2_pt_alpha_zero} on  \pref{ref_eq1}.
Thus, equation \pref{ref_eq1} can't be equivalent to the second
Painlevé equation under point transformations.

{\it The equation \pref{ref_eq2}:}  The specialization of  \pref{p2_pt} on  this equation gives the following transformation~(in~1.11 seconds)
{\small
\begin{system}\nonumber
\bar p &=&  {\Frac {1}{36}}{\Frac {  \left( c+3 \right)  \left( c-2 \right) ^{2}p}
{ \left( 1+c \right)  \left( c-5 \right){y}^{12} }}
\times  (  9\,{y}^{3}c + 66\,{y}^{6}p + \cdots \, 
+27\,{y}^{3} )\bar \alpha\,{\bar y}^{2},\\[5mm]

\bar x &=& {\Frac {2}{3}}{\Frac { \left( -27\,{y}^{6}+3\,{y}^{2}xc-2\,{c}^{2}{p}^{2}-24\,{y}^{6}c+
3\,{y}^{6}{c}^{2}+5\,c{p}^{2}-18\,{y}^{2}x+6\,{p}^{2}-{c}^{3}{p}^{2}+3
\,{y}^{2}x{c}^{2} \right) }{ \left( c-5 \right) {y}^{6}}}{\bar y}^{2},\\[5mm]

{\bar y}^{3} &=&\Frac {1}{16}{\Frac { \left( c-5 \right) }{1+c}}\bar \alpha,\\[5mm]

{\bar \alpha}^{2} &=&
1728 \Frac{\left( 5-c \right)  \left( 1+c \right) ^{2}} {\left( c+3 \right)  \left( c-2 \right) ^{2}} {y}^{18}
\times
( 
-9{y}^{3}c-66{y}^{6}p+54{y}^{2}px+18{y}^{6}{c}^{2}p-48{y}^{6}pc+18{y}^{2}x{c}^{2}p
\\[4mm]
&& +2\,{c}^{3}{p}^{3}-2\,{p}^{3}{c}^{2}+72\,{y}^{2}pxc-34\,{p}^{3}c-30\,{p}^{3}
-27\,{y}^{3}) ^{-2} .
\end{system}%
}%
For the particular values  of $c$ for which $\bar x$ does not depend on $p$ we obtain
  {\sf division by zero errors} when computing the others components and then the equation  \pref{ref_eq2} can not be mapped
  to Painlevé two with $\alpha\neq 0$. However, the specialization of \pref{p2_pt_alpha_zero} on  \pref{ref_eq2} gives 
  {\small
\begin{system}\nonumber
\bar p &=& -\Frac {4}{9}{\Frac { \left( c-2 \right) ^{2}p \left( c+3 \right) 
 \left( 18\,x{y}^{2}p{c}^{2}   -90\,p{y}^{6} \, + \cdots + \,
18\,{y}^{6}p{c}^{2} \right) }{
 \left(c -5\right) ^{2}{y}^{12}}} {\bar y}^{5},\\[5mm]
 
 \bar x &=& {\Frac {2}{3}}{\Frac { \left( -23\,{y}^{6}+3\,{y}^{2}xc-2\,{c}^{2}{p}^{2}-20\,{y}^{6}c+
3\,{y}^{6}{c}^{2}+5\,c{p}^{2}-18\,{y}^{2}x+6\,{p}^{2}-{c}^{3}{p}^{2}+3
\,{y}^{2}x{c}^{2} \right) }{ \left( c-5 \right) {y}^{6}}}{\bar y}^{2},\\[5mm]

{\bar y}^{6} &=& -{\Frac {27}{4}}\Frac{ \left( c-5 \right) ^{3}}{\left( c+3 \right)  \left( c-2 \right) ^{2} }
{y}^{18}\times( 18\,x{y}^{2}p{c}^{2}-90\,p{y}^{6}-9\,{y}^{3}c+54\,px{y}^{2}-27\,{y}^{3}+72\,pcx{y}^{2}-34\,{p}^{3}c\\[5mm]
&&-72\,pc{y}^{6}+2\,{c}^{3}{p}^{3}-30\,{p}^{3}-2\,{p}^{3}{c}^{2}+18\,{y}^{6}p{c}^{2}) ^{-2}
 \end{system}%
 }%
which is point transformation only when $c=-1$. In this case, the resulting transformation~is
$$
\bar p=4\,{\frac {{\bar y}^{5}p}{{y}^{9}}},\, \bar x=2\,{\frac {{\bar y}^{2}x}{{y}^{4}}},\, 
{\bar y}^{6}=1/4\,{y}^{12}
$$
and this   maps the equation \pref{ref_eq2} to Painlevé two $y''=2 y^{3}+  xy$ (with $\alpha=0$).

\end{example}

\section{Acknowledgment}
The author wishes to thank Michel Petitot  for fruitful discussions, Samuel Vidal for encouragements, Pavel Exner for
financial support and anonymous referees  for their valuable comments.

\bibliographystyle{elsart-harv}
\bibliography{/Users/raoufdridi/Scientifique/Latex/c_db/c}

\end{document}